\documentclass[a4paper, 12pt]{article}
\usepackage{enumerate, theorem}
\usepackage{amsmath, amsfonts, amssymb}
\usepackage[height=22.5cm, width=15cm]{geometry}
\usepackage[all]{xy}
\usepackage{mathrsfs, yfonts}
\usepackage[utf8]{inputenc}

\newcommand{\parag}[1]{\paragraph{\sc{#1.}}}

\newtheorem{thm}{Theorem}[subsection]
\newtheorem{defn}[thm]{Definition}
\newtheorem{cor}[thm]{Corollary}
\newtheorem{prop}[thm]{Proposition}
\newtheorem{lemma}[thm]{Lemma}

\begin{document}
\title{Algebraic differential equations of periods integrals.}

\author{Daniel Barlet\footnote{Barlet Daniel, Institut Elie Cartan UMR 7502  \newline
Universit\'e de Lorraine, CNRS, INRIA  et  Institut Universitaire de France, \newline
BP 239 - F - 54506 Vandoeuvre-l\`es-Nancy Cedex.France. \newline
e-mail : daniel.barlet@univ-lorraine.fr}.}

\maketitle

\parag{Abstract} We explain that in the study of the asymptotic expansion at the origin of  a period-integral like $\int_{\gamma_z} \omega/df$ or of a hermitian period
 like $\int_{f=s} \rho.(\omega/df )\wedge( \overline{\omega'/df})$ the computation of the Bernstein polynomial of the "fresco" (filtered differential equation) associated to the
  pair of germs $(f, \omega)$ gives a better control than the computation of the Bernstein polynomial of the full Brieskorn module of the germ of  $f$ at the origin. Moreover, it is easier to compute as it has a better
  functoriality and smaller degree. We illustrate this in the case where $f \in \mathbb{C}[x_0, \dots, x_n]$ has $n+2$ monomials and is not quasi-homogeneous, by giving an explicit simple algorithm
  to produce  a multiple of the Bernstein polynomial when $\omega$ is a monomial holomorphic volume form. Several concrete examples are given.

\parag{AMS Classification} 32 S 25- 32 S 40

\parag{Key words} Period-integral - Hermitian period -  Formal Brieskorn Module - Geometric (a,b)-module - Fresco - Bernstein polynomial.

\newpage

\tableofcontents

\section{Introduction}
This article  simplifies and improves two unpublished papers, see \cite{[B.13]} and \cite{[B.15]}  on the computation of the Bernstein polynomial associated to a  period-integral or to a hermitian period.\\

Let me explain the goal of this paper. We consider a holomorphic function $f : U \to \mathbb{C}$ on an open neighborhood $U$ of the origin in $\mathbb{C}^{n+1}$ which is singular at $0$. The singularity is not assumed to be isolated, but
we choose $U$ small enough in order that $f(0) = 0$ is the only critical value of $f$ on $U$\footnote{recall that for any holomorphic function $f$ singular at the origin, such a $U$  always exists.}. We are interested, for instance,  in the meromorphic extension of the holomorphic functions for $\Re(\lambda) \gg \vert h\vert$:
$$ \frac{1}{\Gamma(\lambda)} \int_U \vert f\vert^{2\lambda}.\bar f^h.\rho.\omega\wedge \bar \omega' $$
where $\omega$ and $\omega'$ are given holomorphic $(n+1)-$forms on $U$, $h$ is in $\mathbb{Z}$ and $\rho$ is a $\mathscr{C}^\infty_c(U)$ function identically $1$ near $0$.\\
We shall assume that for a given integer $q \geq 1$ and for a given $\xi \in \mathbb{Q}$ the meromorphic extension of
\begin{equation*}
  \frac{1}{\Gamma(\lambda)} \int_U \vert f\vert^{2\lambda}.\varphi \tag{$H(\xi, q)$}
  \end{equation*}
have poles at points in $\xi + \mathbb{Z}$ which are of order at most $q-1$, for  any  $(n+1, n+1)$ differential form $\varphi \in \mathscr{C}^\infty_c(U \setminus \{0\})$ (so $\varphi \equiv 0$ near $0$).\\
For instance, if the singularity of $f$ is isolated at the origin, this hypothesis will always be true for $q = 1$ and any $\xi \in \mathbb{Q}$.\\
When the singularity of $f$ is not isolated, this condition will be satisfied for $q = 1$, for instance  as soon as the local monodromy of $f$ acting on the reduced cohomology of the Milnor' fibre at each point nearby $0$ but distinct of $0$ does not present the eigenvalue $\exp(2i\pi.\xi)$.\\
In general, this assumption means that any pole of order $\geq q$ at a point in $\xi + \mathbb{Z}$ comes form the germ of our  situation at the origin. But note that even if the eigenvalue $\exp(2i\pi.\xi)$ of  the local monodromy of $f$ acting on the reduced cohomology of the Milnor' fibre at any point nearby $0$ (including $0$) is simple, we may find a pole of order $2$ for such an integral because of the phenomenon of  "entanglement of consecutive strata" (see \cite{[B.91]} for a topological description and \cite{[B.08]} for a description in term of "Brieskorn modules" of this phenomenon) may appear.\\
We shall give in the theorem \ref{tool 4} and  corollary \ref{fond.3} some necessary numerical conditions which control the order of poles in $\xi + \mathbb{Z}$ \ for a given holomorphic $(n+1)-$form which is much more precise than the "classical condition" (which is in fact valid for such an integral when we replace $\rho.\omega\wedge \omega'$ by any differential form $\varphi \in \mathscr{C}^\infty_c(U)^{n+1, n+1}$) asking that the Bernstein polynomial of $f$ at $0$  has at most $(q-1)$ roots (counting multiplicities) in the set $\xi + \mathbb{Z}$.\\
The precise result for the Mellin transforms of hermitian periods is given in the theorem \ref{tool 4} (the remark following the proof of the corollary \ref{fond.3} indicates also a variant which can be obtained by the same method.)\\
The examples given at the end of this paper show not only that the Bernstein polynomial of the {\bf fresco} associated to the pair $(f, \omega)$  is much easier to compute than the full Bernstein polynomial of $f$, but also that it has, in general, a much smaller number of roots.\\

The main tool around this kind of technic will be the following generalization of the use of a Bernstein identity to control the poles of the Mellin transform of a "hermitian period" of the form 
\begin{equation}
F(\lambda) := \frac{1}{\Gamma(\lambda)}.\int_U \vert f\vert^{2\lambda}.\bar f^h.\omega\wedge \psi 
\end{equation}
where $\omega$ is a  holomorphic form in $\Omega^{n+1}(U)$, where $\psi \in \mathscr{C}^{\infty, (0, n+1)}_c(U)$ is {\bf $d-$closed near $0$} and where $h $ is in $\mathbb{Z}$.

\begin{thm}\label{serieux 0}
Let $f : U \to \mathbb{C}$ be a holomorphic function on a polydisc $U$ with center $0$ in $\mathbb{C}^{n+1}$ and assume that $f(0) = 0$ is the only critical value of $f$ on $U$. For $\xi$ given in $\mathbb{Q}$, assume that the hypothesis $H(\xi, 1)$ is satisfied (see above).\\
 Let $\omega$ be a $(n+1)$ holomorphic differential form on $U$. 
  Assume that the class induced by  $\omega$ in $E^{n+1}$ is annihilated by the element\footnote{Compare this hypothesis with the theorem \ref{frescos.0} which is recalled in section 2.} 
$$P := (a-\lambda_1.b).S_1(b)...(a-\lambda_k.b).S_k(b)$$
in $\tilde{\mathcal{A}}$, where for each $j \in [1, k], \ S_j \in \mathbb{C}[[b]]$ satisfies $S_j(0) = 1$. Note that this assumption depends only of the germs at the origin of $f$ and of $\omega$.\\
Now fix  $\psi \in \mathscr{C}_c^{\infty,(0, n+1)}(U)$ which is  {\bf $d-$closed near $0$}  and assume that for some $h \in \mathbb{Z}$ the meromorphic extension of 
$$ F_h^\psi(\lambda)[\omega] := \frac{1}{\Gamma(\lambda)}.\int_U \vert f\vert^{2\lambda}.\bar f.\omega\wedge \psi $$
has a pole of order $d \geq 1$ at some point in $\xi+ \mathbb{Z}$. Then there exists at least $d$ values of $j \in [1, k]$ such that $\lambda_j$ is in $\xi + \mathbb{Z}$.\\
\end{thm}

Note that the hypothesis of the existence of such a $P$ is always true. But in practice (see the section 4) we may have such a $P$ but we dont know that its initial form in (a,b) corresponds to the Bernstein polynomial of the fresco associated to $\omega$. It is only a member of the (principal)  left ideal which annihilates  the class of $\omega$ in $E^{n+1}$.\\
So, under the hypothesis of this theorem,  the Bernstein polynomial of the fresco $F_\omega$ associated to the pair $(f, \omega)$ (which is the  geometric (a,b)-module generated by $[\omega]$ in $E^{n+1}$)  divides the polynomial 
$$B(\lambda) := \prod_{j=1}^k (\lambda+\lambda_j+j-k).$$

In the more precise statement given in section 3 (see theorem \ref{serieux 1})  we precise the values of these roots of the Bernstein polynomial of $F_\omega$ from the "jumps" of the orders of poles in $\xi + \mathbb{Z}$.\\

For the definition of the geometric (a,b)-module $E^{n+1}$ and  the notion of the Bernstein polynomial of a fresco, see the section 2 below.\\

We shall give a more precise result in the theorem \ref{serieux 1} and some interesting variants using the hypothesis $H(\xi, q)$  in the  theorem \ref{tool 4} and corollary \ref{fond.3} in section 3.

\parag{Remark} Notice that here we use in fact only a "one variable" differential equation (in fact multiplication by the variable in $\mathbb{C}$ and integration relative to this variable) instead of partial differential operators on $\mathbb{C}^{n+1}$ as in the Bernstein identity for $f$ at the origin. This is precisely one of  the interest of using  an (a,b)-module structure in this setting.\\

So we begin this article by  a short overview on geometric (a,b)-modules and frescos intended for the reader which is not familiar with the use of Brieskorn modules in the study of the singularities of
a holomorphic function on a complex manifold.\\
 In opposition with the preprint \cite{[B.15]} cited above, we let aside the global point of view, that is to say the study of the global fresco associated to a period-integral in the case of a proper holomorphic function
on a complex manifold, because it uses more heavy tools and very often the local study presented here would be enough to obtain good informations, using a partition of unity.\\ 
The reader interested in this global setting may consult the  preprint \cite{[B.15]} mentioned above and also the preprint  \cite{[B.12]}.\\
 It is important to notice that we are dealing here with general singularities of a holomorphic function (not only the isolated singularity case  as in the classical use of the Brieskorn module) and 
our illustration in the case of a (not quasi-homogeneous) polynomial in $\mathbb{C}[x_0, \dots, x_n]$ with $n+2$ monomials does not assume also that the singularity is isolated.

\section{A short overview on (a, b)-modules and frescos}

\subsection{Why to use an $(a, b)-$module structure instead of a differential system ?}

Note first that in "$(a, b)$"  $a$ is the multiplication by the variable $z$ and $b$ is the primitive vanishing at $z = 0$, so $b(f)(z) := \int_0^z f(t).dt$ where $f$ is, for instance, a holomorphic
multivalued function with  an eventual ramification point at $z=0$. So we are working with the non commutative algebra $\mathcal{A}$ generated by $a$ and $b$ with the commutation relation
$a.b - b.a = b^2$  as unique relation. This relation corresponds to the usual commutation relation $\partial_z.z-  z.\partial_z = 1$ in the Weyl algebra $\mathbb{C}\langle z, \partial_z\rangle$.
\smallskip

Then why not to use the usual Weyl algebra ?

\smallskip

The initial motivation comes from the study of germs of isolated singularities of holomorphic functions $(f, 0) : (\mathbb{C}^{n+1}, 0) \to (\mathbb{C}, 0)$ initiated at the end of the sixties by Milnor \cite{[Mi.68]}, Brieskorn \cite{[Br.70]}, Deligne \cite{[De.70]},   Malgrange, \cite{[Ma.74]}, Varchenko \cite{[V.80]}, Saito Kyoji  \cite{[Sk.83]}, Saito Morihiko \cite{[Sa.89]}, ... and many others.\\
To my knowledge the first who introduced the "operator" $\partial_z^{-1}$ was Kyoji Saito in the beginning of the eighties (see \cite{[Sk.83]}). The main reason comes from the fact that, looking at period-integrals of the type
$z \mapsto \int_{\gamma_z} \omega/df$ where  $ (f, 0) : (\mathbb{C}^{n+1}, 0) \to (\mathbb{C}, 0)$ is a germ of holomorphic function with an isolated singularity, $\omega \in \Omega^{n+1}_{\mathbb{C}^{n+1}, 0}$ is a germ of holomorphic volume form and $\gamma_z, z \in H$ is  a horizontal family of compact $n-$cycles in the fibers $\{f = z\}$ of $f$, $H \to D^*$ been the universal cover of a small punctured disc $D^*$ with center $0$.\\
The derivation $\partial_z$ of such an integral is given by the following formula
\begin{equation}
 \partial_z\big(\int_{\gamma_z} \omega/df \big) = \int_{\gamma_z} du/df 
 \end{equation}
 where $u \in \Omega^n_{\mathbb{C}^{n+1}, 0}$ satisfies $\omega = df\wedge u$. But in general it is not possible to find such a $u \in \Omega^n_{\mathbb{C}^{n+1}, 0}$ because  writing $\omega = g(x).dx$, the holomorphic germ $g$ is not  in the Jacobian ideal of $f$. Nevertheless, as the coherent sheaf $\Omega^{n+1}\big/df\wedge \Omega^n$ has support inside $\{ f = 0\}$ near $0$, the Nullstellensatz gives a positive integer $p$ such that $f^p$ annihilates this sheaf near $0$ and we may find $u \in f^{-p}.\Omega^n_{\mathbb{C}^{n+1}, 0}$ such that $\omega = df\wedge u$ and $(2)$ holds. But, of course, this implies that in formula $(2)$ the derivation in $z$ needs a denominator which is a power of $z$.\\
 Thanks to the positivity theorem of Malgrange (see \cite{[Ma.74]}) we may write the formula $(2)$ as follows:
 \begin{equation}
 \int_0^z \big(\int_{\gamma_t} du/df\big).dt = \int_{\gamma_z} \omega/df = \int_{\gamma_z} u
 \end{equation}
If  we begin with $\omega := du \in \Omega^{n+1}_{\mathbb{C}^{n+1}, 0}$ we see that formula $(3)$ does not need any denominator in $z$. Moreover the surjectivity of the de Rham differential $d : \Omega^n \to \Omega^{n+1}$ shows that $du$ may be any germ in $\Omega^{n+1}_{\mathbb{C}^{n+1}, 0}$ and we may write
 \begin{equation}
 b\big(\int_{\gamma_z} du/df \big) = \int_{\gamma_z} df\wedge du/df =  \int_{\gamma_z} u 
 \end{equation}
 so that the action of $b$ only needs to solve the equation $du = \omega$, and this  is always possible with $u \in \Omega^n_{\mathbb{C}^{n+1}, 0}$ without introducing a denominator in $f$ (so no denominator  in $z$ downstairs).\\
 As $a$ is given by the formula :
 \begin{equation}
 a\big(\int_{\gamma_z} du/df \big) = z.\big(\int_{\gamma_z} du/df \big) = \int_{\gamma_z} f.du/df \big) 
 \end{equation}
because the $n-$cycle $\gamma_z$ is in the fiber $f^{-1}(z)$ we see that the $\mathcal{A}-$module structure on the quotient\footnote{This quotient allows to define $b[\omega] = [df\wedge u]$ independently on the choice of $u\in \Omega^n_0$ such that $\omega = du$ because when $\omega$ is in $d(Ker(df)_0^n$ the period-integral is identically $0$ near $z = 0$.} $E_f:= \Omega^{n+1}_{\mathbb{C}^{n+1}, 0}\big/ d\big(Ker(df)^n_0\big)$ does not need to consider meromorphic $(n+1)-$differential forms with poles along the fibre  $\{f = 0\}$.\\
 Note that in the case of an isolated singularity for $f$ we have the equality 
 $$Ker(df)^n_0 = df\wedge \Omega^{n-1}_0$$ 
 because the partial derivatives of $f$ define a regular sequence at the origin. Also in this case we find that $E_f\big/bE_f$ is equal to the finite dimensional vector space $\mathcal{O}_{\mathbb{C}^{n+1, 0}}\big/J(f)_0$ where $J(f)$ is the Jacobian ideal of $f$. So we find the classical "Brieskorn lattice".
 
 \smallskip
 
 But why is this presentation interesting if, at the end, we are compelled to introduce denominators in $f$ (or in $z$ working downstairs) to reach an ordinary differential system (or a differential equation ) ?
 
 \smallskip
 
 The answer comes from the following considerations:\\
 If you keep a module structure over the algebra $\mathcal{A}$ as a substitute for a differential system  you have a richer structure (so more precise informations)  than a structure of module over the localized Weyl algebra $\mathbb{C}\langle z, z^{-1}, \partial_z\rangle$ associated to your differential system. This comes from the fact that the commutation relation $a.b - b.a = b^2$ is homogeneous of degree 2 in $(a, b)$ and  implies the existence of the decreasing sequence of two-sided ideals in $\mathcal{A}$ given by $b^m.\mathcal{A} = \mathcal{A}.b^m, \forall m \in \mathbb{N}$. So any $\mathcal{A}-$module $E$ is endowed with a "natural filtration" $(b^m.E)_{m \in \mathbb{N}}$ by sub$-\mathcal{A}-$modules. For instance Varchenkho \cite{[V.80]} proves that in the case of an isolated singularity this filtration defines the Hodge filtration of the mixed Hodge structure on the cohomology of the Milnor's fibre of $f$.
 
 \parag{exercice}  Show that $a.b^m = b^m.a + m.b^{m+1}, \forall m \in \mathbb{N}$ is consequence of the commutation relation corresponding to $m = 1$.\\
 Note that this relation implies  that $a$ and $b^m$ commute modulo $b^{m+1}.\mathcal{A}$.\\

Also,  looking at the "natural action" of $\mathcal{A}$ on $\mathbb{C}[[z]]$ which is given by \  $a(z^m) = z^{m+1} $ and $b(z^m) = z^{m+1}/(m+1)$, you will see that $b^m.\mathbb{C}[[z]] = z^m.\mathbb{C}[[z]], m \in \mathbb{N}$ is the filtration by the valuation.\\

 Another simple remark may also help to convince the reader that a module structure over $\mathcal{A}$ is interesting:
 
 \begin{lemma}\label{construction}
 Let $E := \oplus_{j=1}^k \mathbb{C}[b].e_j$ be a free $\mathbb{C}[b]-$module with basis $e_1, \dots, e_k$ and let $x_1, \dots, x_k$ be any given collection of elements in $E$. Then there exists an unique $\mathcal{A}-$module structure on $E$ such that
 \begin{enumerate}[a)]
 \item The action of $a$ is defined by $a.e_j = x_j$ for each $j \in [1, k]$.
 \item The action of $b$ is given by the $\mathbb{C}[b]-$structure of $E$.
 \end{enumerate}
 \end{lemma}
 
 The proof of this lemma is easily deduced from the following formula which is an easy  consequence of the exercice above:
 $$ a.(S(b).e_j) = S(b).x_j + b^2.S'(b).e_j \quad \forall j \in [1, k]$$
 where $S'(b)$ is the "usual" derivative of the polynomial $S \in \mathbb{C}[b]$. \\
 
 In fact, the presence of the filtration by the two-sided ideals $b^m.\mathcal{A}$ of the algebra $\mathcal{A}$ and the lemma above leads to the following considerations
 \begin{itemize}
 \item The "fundamental" operation\footnote{This is psychologically the most difficult fact to accept after a standard education in maths.} in the action of $\mathcal{A}$ is $b$ !
 \item It seems convenient, as we are interested in the asymptotic expansions of the period-integrals $\int_{\gamma_z} \omega/df$ when $z \to 0$, to complete the algebra $\mathcal{A}$ for the uniform structure defined
 by the filtration $b^m.\mathcal{A},  m \in \mathbb{N}$.\\
 Note that for the "obvious" action of $\mathcal{A}$ on formal power series in $z$ this filtration is associated to the valuation in $z$ (see the remark following the exercise above) .
 \end{itemize}
 This means that we shall work with the algebra
 \begin{equation}
 \tilde{\mathcal{A}} := \{ \sum_{\nu \geq 0} P_\nu(a).b^\nu, \ P_\nu \in \mathbb{C}[a] \quad \forall \nu \in \mathbb{N} \}.
 \end{equation}
 The initial idea of Kyoji Saito was to add some convergence conditions in order that such series acts on convergent (multivalued) series likes
 $$ \sum_{r \in R, j \in [0, N]}\mathbb{C} \{z\}.z^r.(Log \, z)^j $$
 where $R$ is a finite subset in $\mathbb{Q}$ and $N$ is a non negative integer, which are the kind of asymptotic expansions which are valid  for our period-integrals.\\
 But thanks to the regularity of the Gauss-Manin connection, we dont loose any information by staying at the formal series level and this avoids a lot of painful (standard) estimates !\\
  Remark also that the construction given in the lemma above is also valid for the algebra $\tilde{\mathcal{A}}$ and, moreover, that a module $E$  over  $\tilde{\mathcal{A}}$ without $b-$torsion is of  finite type over $\mathbb{C}[[b]]$ if and only if  the complex vector space $E/b.E$ is finite dimensional.\\
  So, our definition of an (a, b)-module will be
 \begin{itemize}
 \item {\bf An (a,b)-module is a left $\tilde{\mathcal{A}}-$module which is a free and finite type module over the (commutative) sub-algebra $\mathbb{C}[[b]] \subset \tilde{\mathcal{A}}$}.
 \end{itemize}
 
 \parag{Examples} 
 \begin{enumerate}
 \item Let $(f, 0) : (\mathbb{C}^{n+1}, 0) \to (\mathbb{C}, 0)$ be a germ of holomorphic function with an isolated singularity. Let $\hat{\Omega}^p_0$ be the formal completion at the origin of $\Omega^p_{\mathbb{C}^{n+1}, 0}$. 
 The quotient $\hat{E}_f := \hat{\Omega}^{n+1}_0\big/\big(df \wedge d(\hat{\Omega}^n_0)\big)$ endowed with the actions of $a := \times f$ and $b := df\wedge d^{-1}$ is an (a,b)-module (Note that the absence of $b-$torsion is a theorem; see \cite{[Se.70]} or \cite{[B.06]}).
\item  Let $E := \mathbb{C}[[b]]e_0 \oplus \mathbb{C}[[b]].e_1$ be the $\tilde{\mathcal{A}}-$module defined by $a.e_0 := b.e_0$ and $a.e_1 = b.e_1- b.e_0$. Then it is an easy exercice to show that $E$ is isomorphic to
$$ \mathbb{C}[[z]] \oplus \mathbb{C}[[z]].Log\, z$$
where $a := \times z$ and $b := \int_0^z $.\\
Determine the filtration $(b^m.E)_{m \in \mathbb{N}}$ in this example and compare it with the filtration by the $(a^m.E)_{m \in \mathbb{N}}$.\\
 Compute the module over the Weyl algebra generated by $Log\, z$ and compare with $E$.\\
\end{enumerate}

\subsection{Geometric (a,b)-modules}

The (a,b)-modules which appear in singularity theory of a function are special. They correspond to {\bf regular} differential system and the notion of regularity is easy to define for an (a,b)-module:\\
First we say that the (a,b)-module $E$ has a {\bf simple pole} when $a.E \subset b.E$. When it is the case, $-b^{-1}.a$ acts on the (finite dimensional) vector space $E/b.E$ and its {\bf minimal polynomial} is called the {\bf Bernstein polynomial} of $E$.\\
For a general (a,b)-module the saturation $\tilde{E}$ of $E$ by the action of $b^{-1}.a$  is not always  a finite type $\mathbb{C}[[b]]-$module. When  $\tilde{E}$ is of finite type over $\mathbb{C}[[b]]$, $\tilde{E}$ is an (a,b)-module (with simple pole) with the same rank over $\mathbb{C}[[b]]$ than the rank of $E$.\\
 We say in this case that $E$ is {\bf regular}. This is equivalent to the fact that $E$ can be embedded in an (a,b)-module having a simple pole.\\
 Then we defined the {\bf Bernstein polynomial} of a regular (a,b)-module $E$  as the Bernstein polynomial of its saturation $\tilde{E}$ by $b^{-1}.a$.\\

There is one more specific property for the (regular) (a,b)-modules coming from the singularity of a function $f$, which reflects the fact that the monodromy of $f$ is quasi-unipotent and the positivity theorem of Malgrange: the fact that the roots of the Bernstein polynomial are negative rational numbers (compare with  the famous theorem of Kashiwara \cite{[K.76]}). So we call {\bf geometric} a regular (a,b)-module whose Bernstein polynomial have negative rational roots.

\parag{Example} In the previous example 2 the (a,b)-module has a simple pole and its Bernstein polynomial is, by definition,  the minimal polynomial of the matrix
 $$\begin{pmatrix}-1 & -1\\ 0 & -1\end{pmatrix}$$
  so its Bernstein  polynomial is $(\lambda + 1)^2$. Compare with the Bernstein type  identity
$$\big((2\lambda+1).\partial_z - z.\partial_z^2\big)(z^{\lambda+1}.Log\, z) = (\lambda+1)^2.z^\lambda.Log\, z .$$

 The following easy proposition will be needed in the sequel. Although it is rather standard, we shall sketch the proof for the convenience of the reader
 
 \begin{prop}\label{geometric}
 Let $E$ be a geometric (a,b)-module and $F$ any sub$-\tilde{\mathcal{A}}-$module in $E$. Then $F$ is a geometric (a,b)-module. When $F$ is normal, the quotient $E/F$ is again a geometric (a,b)-module.
 \end{prop}
 
 \parag{proof} From the regularity of $E$ we may assume that $E$ is a simple pole module (i.e.  $a.E \subset b.E$). Then the Bernstein polynomial of $E$ is the {\bf minimal polynomial} of the action of $-b^{-1}.a$ on the finite dimensional vector space $E/b.E$. As $F$ is a $\mathbb{C}[[b]]$ sub-module of $E$ which is free and finite rank on $\mathbb{C}[[b]]$, $F$ is also free and finite rank on $\mathbb{C}[[b]]$ and stable by $a$. So $F$ is an (a,b)-module. Its saturation by $b^{-1}.a$ is again contained in $E$ and so it is also free  of finite type on $\mathbb{C}[[b]]$. This gives the regularity of $F$. The last points to prove is the fact that the Bernstein polynomial of $F$ has negative rational roots (i.e. $F$ is geometric) and the fact that when $F$ is normal $E/F$ is also geometric. We shall argue by induction on the rank of $F$. In the rank $1$ case let $e$ be a generator of $F$ over $\mathbb{C}[[b]]$ such that $a.e = \lambda.b.e$  (see the classification of rank $1$ regular (a,b)-module in [B.93], lemma 2.4). Let $\nu$ in $\mathbb{N}$ be maximal such that $b^{-\nu}.e$ \ lies in $E$. Then $\mathbb{C}[[b]].b^{-\nu}.e = b^{-\nu}.F$ is a normal\footnote{$G \subset E$ is a normal sub-module when $bx \in G$ for $x \in E$ implies $x \in G$.} sub-module of $E$ and we have an exact sequence of simple poles (a,b)-modules
 $$ 0 \to b^{-\nu}.F \to  E \to  Q \to 0 $$
 and also an exact sequence of $(-b^{-1}.a)$ finite dimensional vector spaces
 $$  0 \to  \mathbb{C}.b^{-\nu}.e \to  E/b.E \to  Q/b.Q  \to 0 .$$
 Then the minimal polynomial $B_{E}$  of the action of $-b^{-1}.a$ on $E/b.E$ is either equal to the minimal polynomial $B_{Q}$ of the action of $-b^{-1}.a$ on $Q/b.Q$, and in this case $-(\lambda -\nu)$ divides $B_{Q} = B_E$,  or we have $B_{E}[x] = (x + (\lambda - \nu)).B_{Q}[x]$. In both cases, as $E$ is geometric, we obtain that $-(\lambda -\nu)$ is a negative rational number, and so is $-\lambda$. Moreover, in both cases, $Q$ is also geometric.\\
 The induction step follows easily by considering a rank 1 normal sub-module $G$ in $F$, using the following lemma and  the fact that a quotient of a geometric (a,b)-module by a normal rank $1$ sub-module is again geometric already proved above. $\hfill \blacksquare$\\
 
 For a proof a the following lemma see for instance the remark 1.2 following the  proposition 1.3 in  \cite{[B.93]}.
 
 \begin{lemma}\label{racine}
Let $E$ be a regular (a,b)-module and let $F \subset E$ be a sub-(a,b)-module. Assume that $\lambda$ is a root of the Bernstein polynomial $B_F$ of $F$. Then there exists $\lambda' \in \lambda + \mathbb{N}$ such that $\lambda'$ is a root of the Bernstein polynomial $B_E$ of $E$. 
\end{lemma}

As we want to consider, in the non isolated singularity case, a sheaf of geometric (a,b)-module along the singular set $\{df = 0 \}$ of the zero set $\{f = 0 \}$ of a holomorphic function on a complex manifold $M$, we have to replace the completion used  in the classical case of an isolated singularity case by a $f-$completion which is in fact the $z-$completion downstairs.  This will not change seriously the considerations above, thanks to the following easy proposition which implies that any geometric (a,b)-module is in fact a module over the algebra $\hat{\mathcal{A}} := \{ \sum_{p, q \geq 0} c_{p, q}.a^p.b^q \}$ which contains both $\mathbb{C}[[b]]$ and $\mathbb{C}[[a]]$.

\begin{prop}\label{a-complete}
Any geometric (a,b)-module is complete for the decreasing filtration by the $\mathbb{C}[a]-$sub-modules (not stable by $b$ in general)  $(a^m.E)_{m \in \mathbb{N}}$.
\end{prop}

Note that the hypothesis "geometric" insure that $a$ is injective (this is not the case if we assume only the regularity). But the condition that there exists $k \in \mathbb{N}^*$ such that $a^k.E \subset b.E$ is enough to see that $\mathbb{C}[[a]]$ acts on the (a,b)-module $E$. This condition is clearly satisfied by any regular (a,b)-module (see \cite{[B.93]}).\\

 \subsection{Frescos}
 
 We have seen that the (a,b)-module structure may be an interesting way to study the differential system associated to period-integrals for a germ of holomorphic function $(f,0) : (\mathbb{C}^{n+1}, 0) \to (\mathbb{C}, 0)$. In fact, the Brieskorn module, or the (a,b)-module $\hat{E}_f$ 
  which is the completion of the Brieskorn module $\Omega^{n+1}_0\big/df\wedge d\Omega^{n-1}_0$ in the isolated singularity case. It gives in fact a filtered version of the differential system satisfied by {\bf all the period-integrals} associated to the germ  $(f, 0)$.\\
 But if we are interested by the period-integrals corresponding to a specific holomorphic differential form, it is clear that such a differential system, that is to say the all (a,b)-module $\hat{E}_f$, does not give very precise informations. In term of differential system, we would prefer to have a specific {\bf differential equation} satisfied by the integral-periods $\int_{\gamma_z} \omega/df$ for our choice of $\omega$ than the differential system satisfied by all period-integrals, so associated to all choices of $\omega \in \Omega^{n+1}_0$. The analog of the differential equation in term of (a,b)-modules is the notion of "fresco". A {\bf fresco} is a, by definition, a {\bf geometric} (a,b)-module which is generated, as a $\tilde{\mathcal{A}}-$module, by {\bf one generator}. For instance, in the previous situation, we shall consider the fresco given by $\tilde{\mathcal{A}}.[\omega] \subset \hat{E}_f$ and we shall call it the {\bf fresco of the pair $(f, \omega)$ at the origin}. The following structure  theorem  describes in a very simple way such a fresco (see \cite{[B.09b]}). 
 
  \begin{thm}\label{frescos.0}
 Any rank $k$ fresco\footnote{Recall that we consider the rank as a $\mathbb{C}[[b]]-$module, where  $\mathbb{C}[[b]] \subset \tilde{\mathcal{A}}$.} $F$ with generator $e$ is isomorphic (as an $ \tilde{\mathcal{A}}-$module) to a quotient $ \tilde{\mathcal{A}}\big/ \tilde{\mathcal{A}}.\Pi$, the isomorphism $F \to  \tilde{\mathcal{A}}\big/ \tilde{\mathcal{A}}.\Pi$ being defined by sending  the generator $e$ of $F$  to the class of  $1$. We may choose $\Pi$  having  the following form
 \begin{equation}
 \Pi := (a - \lambda_{1}.b).S_{1}.(a - \lambda_{2}.b).S_{2} \dots S_{k-1}.(a - \lambda_{k}.b).S_k 
 \end{equation}
 where the numbers  $-(\lambda_{j}+j-k)$ are the roots of the Bernstein polynomial of $F$ and where $S_{j}$ are in $\mathbb{C}[b]$ and satisfy $S_{j}(0) = 1$ (so each $S_j$ is invertible in $\mathbb{C}[[b]]$).
 \end{thm}
 
 Note that the initial form in (a,b) of $\Pi$ is $P_{F} := (a - \lambda_{1}.b)\dots (a - \lambda_{k}.b)$. It is called  the {\bf Bernstein element} of  the fresco $F$. It does not depend of the choice of the generator of $F$ over $\tilde{\mathcal{A}}$ (choice which determines $\Pi$)   and is related to the {\bf Bernstein polynomial} $B_F \in \mathbb{C}[\lambda]$  of $F$ by the relation  in the ring $\mathcal{A}[b^{-1}]$:
 \begin{equation}
  (-b)^{k}.B_{F}(-b^{-1}.a) = P_{F}, \quad{\rm where} \quad k:= rk(F).
  \end{equation}

  In the case of a fresco $F$ the Bernstein polynomial $B_F$ is equal to the {\bf characteristic polynomial} of the action of $-b^{-1}.a$ on $\tilde{F}/b.\tilde{F}$ where $\tilde{F}$ is the saturation of $F$ by $b^{-1}.a$ (see \cite{[B.09b]} theorem 3.2.1). This makes the computation
  of the Bernstein polynomial of a fresco easier  than for a general geometric (a,b)-module, for instance by the use of  the following remark:\\
  
  If $0 \to F \to G \to H \to 0$ is an exact sequence of frescos we have the relation $P_{G} = P_{F}.P_{H}$ (product in $\mathcal{A}$) between the Bernstein elements and this gives  the relation (see \cite{[B.09b]} proposition 3.4.4):
   $$B_{G}(x) = B_{F}(x +rk(H)).B_{H}(x) $$
   between the Bernstein polynomials.\\
   
   Recall that any element $P$ in the algebra $\mathcal{A}$ which is homogeneous of degree $k$ in (a,b) and monic in $a$ may be written $P = (a-r_1.b)\dots(a-r_k.b)$ where $r_1, \dots, r_k$ are complex numbers.
   This equality is not unique but the sequence $r_j+j-k, j\in [1, k]$ depends only on $P$ (see \cite{[B.09b]} proposition 2.0.2).
   
   Our next proposition will be useful in computations of examples.

  \begin{prop}\label{egalite}
  Let $F$ be a rank $k$ fresco with generator $e$. Assume that $\mathcal{Q} \in \hat{\mathcal{A}}$ has the following properties:
  \begin{enumerate} [i)]
  \item The initial form $Q$ of $\mathcal{Q}$ in $(a, b)$ has degree $d$ and is monic in $a$.
  \item $\mathcal{Q}[e] = 0$ in $F$.
  \end{enumerate}
  Then $Q$ is a left multiple in $\mathcal{A}$ of $P_F$, the Bernstein element of $F$.\\
 If moreover we have $d = k$, then  $Q$ is the Bernstein element of $F$.
  \end{prop}
  
  \parag{Proof} Using the structure theorem of \cite{[B.09b]} recalled in the theorem \ref{frescos.0} above, we have an isomorphism  $F \simeq \tilde{\mathcal{A}}\big/\tilde{\mathcal{A}}.\Pi$ where the initial form in (a,b)  $P_F$ of $\Pi$ is the Bernstein element of $F$. As $F$ is a $\hat{\mathcal{A}}-$module (see the proposition \ref{a-complete}) we have also an isomorphism  $F \simeq \hat{\mathcal{A}}\big/\hat{\mathcal{A}}.\Pi$ of $\hat{\mathcal{A}}-$modules and our hypothesis $ii)$ implies that there exists $Z \in \hat{\mathcal{A}}$ such that
  $$ \mathcal{Q} = Z.\Pi .$$
  This gives $Q = in(Z).P$ where $in(Z)$ is the initial form in (a,b) of $Z$. This already implies that $d \geq k$ and that $in(Z)$ is of degree $d-k$. In the case $d = k$ we have $In(Z) = 1$ and $Q = P$. $\hfill \blacksquare$\\

 \subsection{A general existence theorem}
 
 Now consider a germ $f : (\mathbb{C}^{n+1},0) \to (\mathbb{C},0)$ such that $\{ f = 0 \}$ is reduced. 
Let $\hat{\Omega}^{\bullet}$ the formal $f-$completion of the sheaf of holomorphic differential forms on $(\mathbb{C}^{n+1}, 0)$ and let $\hat{K}er\, df^{\bullet}$ be the kernel of the map
 $$ \wedge df : \hat{\Omega}^{\bullet} \longrightarrow \hat{\Omega}^{\bullet+1}.$$
 Then for any $p \geq 2$ \ the $p-$th  cohomology sheaf\footnote{For $p=1$ we have to replace $\hat{K}er \,df^1$ by a quotient ; see \cite{[B.08]}.}  of the complex $(\hat{K}er\, df^{\bullet}, d^{\bullet})$ has a natural structure of  left $\hat{\mathcal{A}}-$module, where the action of $a$ is given by multiplication by $f$ and the action of $b$ is (locally) given by $df \wedge d^{-1}$. \\
   
  The following result is  known (see \cite{[B.06]}, \cite{[B-S.07]} and \cite{[B.08]})
 
 \begin{thm}\label{Geometric  modules}
 For each integer $p$ the germ at $0$ of the $p-$th cohomology sheaf   of the complex $(\hat{K}er\, df^{\bullet}, d^{\bullet})$ (modified for $p=1$, see the footnote), denoted by $\mathcal{E}^{p}$,  satisfies the following properties:
 \begin{enumerate}[i)]
 \item We have in $\mathcal{E}^{p}$ the commutation relation $a.b - b.a = b^{2}$.
 \item $\mathcal{E}^{p}$ is  b-separated and b-complete (so also a-complete). Then it is a $\tilde{\mathcal{A}}-$module (and also a $\hat{\mathcal{A}}-$module).
 \item There exists an integer $m \geq 1$ such that $a^{m}.\mathcal{E}^{p} \subset b.\mathcal{E}^{p}$.
 \item We have $B(\mathcal{E}^{p}) = A(\mathcal{E}^{p}) = \tilde{A}(\mathcal{E}^{p})$ and there exists an integer $N \geq 1$ such that $a^{N}. A(\mathcal{E}^{p}) = 0$ and \ $b^{2N}.B(\mathcal{E}^{p}) = 0$.
 \item The quotient $E^p := \mathcal{E}^{p}\big/B(\mathcal{E}^{p})$ is a geometric (a,b)-module.
 \end{enumerate}
 \end{thm}
 
 Recall that $B(\mathcal{E})$ is the b-torsion in $\mathcal{E}$, $\tilde{A}(\mathcal{E})$ the a-torsion of $\mathcal{E}$ and $A(\mathcal{E})$ the $\mathbb{C}[b]-$module generated by $\tilde{A}(\mathcal{E})$ in $\mathcal{E}$. \\
 
 We shall mainly use this result in the case where $\omega$ is an $(n+1)-$holomorphic differential form in an open neighborhood $U$ of the origin of $0 \in \mathbb{C}^{n+1}$; note that the condition $df\wedge \omega = d\omega = 0$ is automatic in this case. So it defines a class $[\omega]$ in $E^{n+1}$ and generates the fresco
  $$F_{\omega} := \tilde{\mathcal{A}}.[\omega] \subset E^{n+1}_0$$
  thanks to the proposition \ref{geometric} and property v) of the previous theorem.

 \begin{defn}\label{fresco omega}
 We shall denote $B_{\omega}\in \mathbb{C}[x]$ and $P_{\omega}\in \mathcal{A}$ respectively the Bernstein polynomial and the Bernstein element of the fresco $F_{\omega}$.
 \end{defn}
 
We shall study several examples   in section 4.

\section{Mellin transform of hermitian periods}

\subsection{The main result}

We consider now a holomorphic function on an open polydisc  $U$  centered at the origin in $\mathbb{C}^{n+1}$ such that $f(0) = 0$ is the only critical value of $f$ on $ U$.
Let $\omega$ be a holomorphic $(n+1)-$differential form on $ U$ . Then let $\psi$ be a $\mathscr{C}^\infty$ differential form with compact support in the  polydisc $ U$ of type $(0, n+1)$ which satisfies $d\psi \equiv 0$ near $ 0$. For any $h \in \mathbb{Z}$ define, at least for $\Re(\lambda)$ large enough, the holomorphic function
\begin{equation}
\lambda \mapsto  F_h^\psi(\lambda)[\omega] = \frac{1}{\Gamma(\lambda)}.\int_U \vert f\vert^{2\lambda}.\bar f^h.\omega\wedge \psi .
\end{equation}
Note that the existence of a Bernstein identity for $f$ in a neighborhood of $0$ insures that for $U$ small enough and any differential form  $\varphi \in \mathscr{C}^\infty_c(U)^{(n+1, n+1)}$ the holomorphic function defined by $\int_U \vert f\vert^{2\lambda}.\varphi$ has a meromorphic extension to the all complex plane with a finite series of poles of order at most $n+1$ at points of the form $\xi_j - \mathbb{N}$ where $\xi_j$ are negative rational numbers (see \cite{[K.76]}). \\

We shall make the following hypothesis, which we shall denote by $H(\xi,1)$:
\begin{itemize}
\item For $\xi \in \mathbb{Q}$ the local monodromy of $f$ acting on the reduced cohomology of the Milnor's fibre at any point distinct from $0$ does not admit the eigenvalue $\exp(2i\pi.\xi)$.
\end{itemize}

\begin{prop}\label{tool 1} 
We assume the hypothesis $H(\xi, 1)$. Let $\omega$ be a holomorphic $(n+1)-$differential form on $ U$  and fix $\psi$ be a $\mathscr{C}^\infty$ differential form with compact support in  $ U$ of type $(0, n+1)$ which satisfies $d\psi = 0$ near $0$.  We have the following formulas in $E^{n+1}$ (recall that the geometric (a,b)-module $E^{n+1}$ is defined in the theorem \ref{Geometric  modules} above and  that the action of $a$ and $b$ are defined by:  $a.[\omega] = [f.\omega]$, $b.[\omega] = [df\wedge u]$ where $u \in \Omega^p(U)$ satisfies $du = \omega$\footnote{Such a $u$  always exists.}).\\
\begin{enumerate}[i)]
\item If there exists  $v \in \Omega^n(U)$ satisfying $df \wedge v \equiv 0$ and $dv = \omega$ on $U$, then $F^\psi_h[\omega] $ has no pole in $\xi + \mathbb{Z}$ for any $h$ and any $\psi$.
\item $F^\psi_h(\lambda)[a.\omega] - (\lambda+1).F^\psi_{h-1}(\lambda+1)[\omega]$ has no pole in $\xi + \mathbb{Z}$.
\item $F^\psi_h(\lambda)[b.\omega]  + F^\psi_{h-1}(\lambda+1)[\omega] $ has no pole in $\xi + \mathbb{Z}$.
\item So for any $\mu \in \mathbb{C}$ the meromorphic function 
 $$F^\psi_h(\lambda)[(a - \mu.b).\omega] - (\lambda+\mu+1).F^\psi_{h-1}(\lambda+1)[\omega]$$
 has no pole in $\xi + \mathbb{Z}$, combining $ii)$ and $iii)$.
\end{enumerate}
\end{prop}

Of course, the simplest example of such a $\psi$ is  given by $\psi := \rho.\bar\omega'$ where $\rho$ is a $\mathscr{C}^\infty$ function with compact support such that $\rho \equiv 1$ in a neighborhood of $0$ and where $\omega'$ is a holomorphic $(n+1)-$differential form in $ U$.

\parag{proof} Write $\omega = du$ on $ U$ with $u \in \Omega^n(U)$. This is always possible as $U$ is a polydisc so Stein and contractible. Then for $\Re(\lambda) > 1+ \vert h\vert $ the differential form
$\alpha := \vert f\vert^{2\lambda}.\bar f^h.u\wedge \psi$ is $\mathscr{C}^1$ in $ U$ and has compact support. So  as we have:
 $$d\alpha = \vert f\vert^{2\lambda}.\bar f^h.du\wedge \psi + \lambda.\vert f\vert^{2(\lambda-1)}.\bar f^{h+1}.df\wedge u\wedge \psi + (-1)^n. \vert f\vert^{2\lambda}.\bar f^h.u\wedge d\psi. $$
 Stokes formula gives, as $d\psi$ vanishes near $0$, 
$$\frac{1}{\Gamma(\lambda)}.\int_U d\alpha = 0 = F^\psi_h(\lambda)[\omega] + F^\psi_{h+1}(\lambda-1)[b.\omega] + G(\lambda)$$
where $G(\lambda)$ is a meromorphic function on $\mathbb{C}$ which has no pole in $\xi + \mathbb{Z}$ thanks to our hypothesis $H(\xi,1)$.
This implies $i)$ and $iii)$ using $\Gamma(\lambda) = \lambda.\Gamma(\lambda-1)$. The formula $ii)$ is easy and left to the reader.$\hfill \blacksquare$\\

\parag{Remark} The point $i)$ of the previous proposition shows that $F_h^\psi(\lambda)[\omega]$ has no pole in $\xi + \mathbb{Z}$ when $\omega$ induces the zero class in $\mathcal{E}^{n+1} = H^{n+1}(\hat{K}er\, df^\bullet_0, d^\bullet)$, and the point $iii)$ implies the same conclusion when the class of $\omega$ in $\mathcal{E}^{n+1} $ is of $b-$torsion. So the polar part of $F_h^\psi(\lambda)[\omega]$ in $\xi + \mathbb{Z}$ depends only of the class induced by $\omega$ in $E^{n+1} = \mathcal{E}^{n+1}\big/b-torsion$ \ for $h$ and $\psi$ fixed. This remark will be crucial in the sequel.\\

We shall give the following more precise result than the theorem \ref{serieux 0} stated in the introduction. Let me recall the situation.
Let $f $ be a holomorphic function on a   polydisc $ U$ with center $0$ in $\mathbb{C}^{n+1}$ and assume that $f(0) = 0$ is the only critical value of $f$ on $ U$. Let $\omega$ be a $(n+1)$ holomorphic differential form on   $ U$ and let $\psi$ be a $\mathscr{C}^\infty$ differential form with compact support in $ U$ of type $(0, n+1)$ which satisfies $d\psi \equiv 0$ near $0$.  Then define as above
\begin{equation*}
F_h^\psi(\lambda)[\omega] := \frac{1}{\Gamma(\lambda)}.\int_U \vert f\vert^{2\lambda}.\bar f^h.\omega\wedge \psi  \tag{9}
\end{equation*}

\begin{thm}\label{serieux 1}
Fix $\xi \in \mathbb{Q}$ and an integer $d \geq 1$. Assume that the hypothesis $H(\xi, 1)$ is satisfied. Then assume that, for a given $\psi$, the meromorphic extension of $F_h^\psi(\lambda)[\omega]$ has no pole of order $\geq d+1$ at any point in  $\xi + \mathbb{Z}$, for any choice of $h \in \mathbb{Z}$, but that there exists a point in $\xi + \mathbb{Z}$
where, for some $h$,  this meromorphic extension has a pole order $d$.\\
 For each integer $s \in [1, d]$ let $\xi_s$ be the biggest element in $\xi + \mathbb{Z}$ for which there exists $h \in \mathbb{Z}$ such that $ F_h^\psi(\lambda)[\omega]$ has a pole of order at least equal to $s$ at $\xi_s$. Then  each  $\xi_s$ for $s \in [1, d]$ is a  root of the Bernstein polynomial of the the fresco $F_\omega$.\\
  Moreover,  if $\xi_s = \xi_{s+1} = \dots = \xi_{s+p}$ then $\xi_s$ is  a  root of the Bernstein polynomial of $F_\omega$ with multiplicity at least equal to $p+1$. 
 \end{thm}
 
 Remark that  the theorem implies that the Bernstein polynomial of the fresco $F_\omega$ is a multiple of $\prod_{s=1}^d (\lambda - \xi_s)$ and that we have $\xi_d \leq \xi_{d-1} \leq \dots \leq \xi_1 < 0$ by definition.\\
 
 The proof of this theorem needs some lemmata.

\begin{lemma}\label{Tool 1}
In the situation of the theorem \ref{serieux 1}, let $S \in \mathbb{C}[[b]]$ which satisfies  $S(0) = 1$ and let $\mu \in \mathbb{C}$ such that $\mu \not= -\xi_s$ for a given $s \in [1, d]$. Then
\begin{enumerate}
\item $\xi_s$ is still the biggest pole of order $\geq s$ in $\xi + \mathbb{Z}$  for the meromorphic extension of $F_h^\psi(\lambda)[S(b)\omega]$ when $h$ varies in $\mathbb{Z}$.
\item $\xi_s -1$ is the biggest pole of order $\geq s$ in $\xi + \mathbb{Z}$ for the polar part of the meromorphic extension of $F_h^\psi(\lambda)[(a-\mu.b)\omega]$ when $h$ varies in $\mathbb{Z}$.
\end{enumerate}
\end{lemma}

\parag{proof} Write $S(b) := 1 + \sum_{m = 1}^\infty  s_m.b^m$. As $F_h^\psi(\lambda)[b^m.\omega] = (-1)^m.F_{h-m}^\psi(\lambda+m)[\omega]$, for each $m \geq 1$, thanks to point $iii)$ in the proposition \ref{tool 1},  the meromorphic extension of $F_h^\psi(\lambda)[b^m.\omega] $ cannot have a pole of order $\geq s$ at the point $\xi_s$ for any choice of $h$. So the pole of order $\geq s$  given by the initial term (i.e. $m = 0$)  for a suitable value of $h$, stays maximal. This proves 1.\\
Because we assume $\mu + \xi_s \not= 0$, the point $iv)$ in the proposition \ref{tool 1} shows that the pole at $\xi_s-1$ of $F_h^\psi(\lambda)[(a-\mu.b)\omega]$ has the same order than the pole at $\xi_s$ for $F_{h+1}^\psi(\lambda)[\omega] $.  Also, the same formula shows that for any integer $p \geq 0$ the order of the pole $\xi_s+p$ for  $F_h^\psi(\lambda)[(a-\mu.b)\omega]$ is less or  equal  to the order of the pole  $\xi_s +p+1 \geq \xi_s +1$ for  $F_{h-1}^\psi(\lambda)[\omega]$ which is at most $s-1$ by definition of $\xi_s$. This allows to conclude.$\hfill \blacksquare$\\

\begin{lemma}\label{Tool 2}
In the situation of the theorem \ref{serieux 1}, assume that $\xi_{s+1} = \xi_s$, for some $s \in [1, d-1]$. Then $\xi_s-1$ is  the biggest pole of order $\geq s$  in $\xi + \mathbb{Z}$  for the meromorphic extension of $F_h^\psi(\lambda)[(a+\xi_s.b)\omega]$ when $h$ varies in $\mathbb{Z}$. 
\end{lemma}

\parag{proof} Using the point $iv)$ of  the proposition \ref{tool 1}  we obtain that, for some suitable choice of $h$, the meromorphic extension of $F_h^\psi(\lambda)[(a+\xi_{s+1}.b).\omega]$ has a pole of order $\geq s$ at the point $\xi_{s+1} - 1 = \xi_s - 1$. Assume that for some integer $p \geq 0$ and some $h \in \mathbb{Z}$ the meromorphic extension of $F_h^\psi(\lambda)[(a+\xi_{s+1}.b).\omega]$ has a pole of order $\geq s$ at $\xi_s + p$. Then using again the formula of the point $iv)$ in  the proposition \ref{tool 1} we find that $(\lambda -\xi_s +1).F_{h-1}^\psi(\lambda+1)[\omega]$ has, for a suitable choice of $h \in \mathbb{Z}$,  a pole of order $\geq s$ at $\lambda = \xi_s + p$. But  $\xi_s + p -\xi_s + 1 = p+1 \not= 0$ so we find a pole of order $\geq s$ at the point $\xi_s + p + 1$ for $F_{h-1}(\lambda)[\omega]$. As $p+1 \geq 1$ this contradicts the definition of $\xi_s$, so $\xi_s - 1$ is the biggest pole of order $\geq s$  in $\xi + \mathbb{Z}$  when $h$ varies in $\mathbb{Z}$, for the meromorphic extension of $F_h^\psi(\lambda)[(a+\xi_s.b).\omega]$.$\hfill \blacksquare$\\

\begin{lemma}\label{Tool 1bis}
In the situation of the theorem \ref{serieux 1}, there exists a minimal integer $j_d \in [0, k-1]$ such that $\xi_d -j_d =  -\lambda_{k-j_d}$ and 
for each $s \in [1, d-1]$ there exists a minimal integer $j_s \in [j_{s+1}+1, k-1]$ such that $\xi_s -j_s = -\lambda_{k-j_s}$. Moreover the meromorphic extension of $F_h^\psi(\lambda)[P_{k-j_s-1}.\omega]$ has a pole of order $\geq s-1$ at the point $\xi_s -j_s -1$, where we define
$$P_r := (a-\lambda_1.b).S_1(b) \dots S_{r-1}.(a-\lambda_r.b).S_r \quad {\rm for} \quad  r \in [1, k]$$
where $P_k := \Pi$ is given by applying the theorem 2.3.1 to the fresco $F_\omega$ and  satisfies moreover the condition that the sequence $(\lambda_j+j), j \in [1, k]$ is non decreasing.
\end{lemma}

\parag{proof} The fact that we may assume that the sequence $(\lambda_j+j), j \in [1, k]$ is non decreasing  is consequence of the existence of a {\bf principal Jordan-H\"older sequence} for a fresco (see \cite{[B.13a]} theorem 1.2.5 or \cite{[B.09b]} corollary 3.5.4).
Note also  that the roots of the Bernstein polynomial of $F_\omega$ are given by the numbers $-(\lambda_j + j - k)$ for $j \in [1, k]$ .\\
 Now the first point is to prove that there exists an integer $j \in [0, k-1]$ such that $\xi_d -j = -\lambda_{k-j}$, because we may then define $j_d$ as the minimal such integer. So assume that no such $j$ exists. Then applying the lemma \ref{Tool 1} we will obtain that  the meromorphic extension of $F_h^\psi(\lambda)(\Pi.\omega)$ has a pole of order at least equal to $d \geq 1$ at the point $\xi_d - k$. But points $i)$ and $iii)$ in the proposition \ref{tool 1} implies that there is no pole in $\xi + \mathbb{Z}$ in the meromorphic extension of $F_h^\psi(\lambda)[\Pi.\omega]$ as $[\Pi.\omega] = 0$ in $E^{n+1}$. Contradiction. So such a $j$ exists and $j_d$ is well defined. \\
The same argument than for $s = d$ shows that for each $s \in [1, d-1]$ there exists at least one $j \in [1, k]$ such that $\xi_s -j = -\lambda_{k-j}$. Now we have
$$ \xi_{s+1} = -\lambda_{k-j_{s+1}} + j_{s+1} \leq \xi_s = -\lambda_{k-j} + j .$$
The non decreasing property of the sequence $\lambda_j + j$ implies then that $j \geq j_{s+1}$. If the inequality is strict, then we obtain  $j_s := j$ and the point $iv)$ in the proposition \ref{tool 1} implies that $F_h^\psi(\lambda)[P_{k-j_s-1}.\omega]$ has a pole of order $\geq s-1$ at the point $\xi_s -j_s -1$.\\
If we have $\xi_{s+1} = \xi_s$ then the lemma \ref{Tool 2} shows that $\xi_s - j_{s+1} -1$ is still in this case  the biggest pole of order $\geq s$ in $\xi + \mathbb{Z}$ for the meromorphic extension of $F_h^\psi(\lambda)[(a+\xi_s.b).P_{k-j_{s+1} -1}.\omega]$ when $h$ varies in $\mathbb{Z}$. So we can continue to apply $S_{k-j_{s+1}-2}$ and then $a - \lambda_{k-j_{s+1}-2}.b, \ etc...$  until we reach another $j < j_{s+1}$ such that $\xi_s -j = -\lambda_j$, and then we conclude using the same argument as above.$\hfill \blacksquare$\\

\parag{Remarks}
\begin{enumerate}
\item The sequence $j_s, s \in [1, d]$ is strictly decreasing, so the sequence $\xi_s -j_s$ is strictly decreasing and there are exactly $d$ rational numbers 
$$\xi_s = -(\lambda_{k-j_s} -j_s) = -(\lambda_{k-j_s} + k-j_s - k)$$
 counting multiplicities (remark  that the multiplicities correspond to equalities $\lambda_{k-j_s} - \lambda_{k-j_{s+1}} = j_s - j_{s+1}$).
\item As long as $r \leq j_s$ the rational number $\xi_s$ stays the biggest pole of order $\geq s$ in $\xi + \mathbb{Z}$ for the meromorphic extension of $F_h^\psi(\lambda)[Q_r.\omega]$ when $h$ describes $\mathbb{Z}$, where 
\begin{equation*}
Q_r :=  S_{r}.(a-\lambda_{r +1}.b).S_{r+1} ... (a-\lambda_k.b).S_k \tag{@}
\end{equation*}
\item Note that if we have several $\psi$ which are $d-$closed near $0$ and for which the meromorphic extension of $F^\psi_h[\omega]$ presents poles in $\xi + \mathbb{Z}$, we may obtain more roots in $\xi + \mathbb{Z}$ for the Bernstein polynomial of $F_\omega$ from this result.
\end{enumerate}

\parag{Proof of the theorem \ref{serieux 1}} We shall prove first,  by induction on $d\geq 1$, that there exists at least $d$ values of $j \in [1, k]$ such that $-\lambda_j$ belongs to $\xi + \mathbb{Z}$.\\
 So  assume that  either $d = 1$ or that $d \geq 2$  and that our claim is proved for $d-1$. Then consider the poles of the meromorphic extension of  $F_h(\lambda)[Q_{k-j_d}.\omega]$  (where $Q_r$ is defined in the formula $(@)$ above) and where  the integer $j_d$ is defined in the lemma \ref{Tool 1bis}. Using the lemma \ref{Tool 1} applied to $\xi_d$, we obtain that it has a maximal pole of order $d$ at the point $\xi_d- j_d$ for a suitable choice of $h$ and, applying the lemma \ref{Tool 1bis} we conclude that the meromorphic extension of  $F_h(\lambda)[(a-\lambda_{k-j_d}.b).Q_{k-j_d}.\omega]$ has a pole of order at least equal to $d-1$ at the point $\xi_d -j_d -1$. But the form $\omega' := (a-\lambda_{k-j_d}.b).Q_{k-j_d}.\omega$ is killed in $E^{n+1}$ by $P_{k-j_d-1}$, so the induction hypothesis gives at least $d-1$ values of $j \in [1, k-j_d -1]$ such that $-\lambda_j$ is in $\xi + \mathbb{Z}$. As $-\lambda_{k-j_d} = \xi_d - j_d $ belongs to $\xi + \mathbb{Z}$ this completes the proof of our induction.\\
So we obtain that at least $d$ roots (counting multiplicities) of the Bernstein polynomial of $F_\omega$ are among the $\xi_s$ for $s \in [1, d]$.$\hfill \blacksquare$\\

\subsection{ Some variants}

The following variant of the previous result is also useful.\\

Assume now that we fix $\xi +\mathbb{Z}$, and integer $q \geq 1$  and that we make the following hypothesis:\\
\begin{itemize}
\item For  any differential form $\varphi \in \mathscr{C}^{\infty, (n+1, n+1)}_c(U \setminus \{0\})$ the meromorphic extension of 
\begin{equation*}
\frac{1}{\Gamma(\lambda)}.\int_U \vert f \vert^{2\lambda}.\varphi \tag{$H(\xi, q)$}
\end{equation*}
has never a pole of order at least equal to $q$ at some point in $\xi + \mathbb{Z}$.\\
Note that we require that  $\varphi \equiv 0$ near $0$.
\end{itemize}

\begin{thm}\label{tool 4} Let $\omega \in \Omega^{n+1}(U)$ and fix a differential form   $\psi \in \mathscr{C}_c^{\infty, (0, n+1)}(U)$ which is { \bf $d-$closed near $0$}.\\
Assume that the condition  $H(\xi, q)$ is satisfied and that for some integer $h \in \mathbb{Z}$ the meromorphic extension of 
\begin{equation}
 F^\psi_h(\lambda[\omega]) := \frac{1}{\Gamma(\lambda)}.\int_U \vert f \vert^{2\lambda}.\bar f^h.\omega\wedge \psi 
 \end{equation}
 has a pole of order $q+d-1$ at some point in $\xi + \mathbb{Z}$ with $d\geq 1$, maximal in $\mathbb{N}^*$. Let $\xi_0$ be maximal in $\xi + \mathbb{Z}$ such there exists $h \in \mathbb{Z}$ with the property that $ F^\psi_h(\lambda)[\omega]$ has a pole of order $ q+d-1$ at $\xi_0$. Then $\xi_0$ is a root of  the Bernstein polynomial $B_\omega$ of the fresco $F_{\omega}$ associated to the germs of $f$ and $\omega$ at the origin, and $B_\omega$ admits at least $d$ roots in $\xi_0 + \mathbb{N}$ (counting multiplicities).
\end{thm}

\parag{proof} The argument is analog to the one given in the proof of the theorem \ref{serieux 1}. The  change to make in the proof is that we must take in account here only the polar parts of order at least equal to $q$  of the poles at points in $\xi + \mathbb{Z}$.  So in the proposition \ref{tool 1} and in the lemmata \ref{Tool 1}, \ref{Tool 2} and \ref{Tool 1bis} we have to replace "no poles in $\xi + \mathbb{Z}$" by  "no pole of order $\geq q$ in $\xi + \mathbb{Z}$" under the hypothesis $H(\xi, q)$. Also we define $\xi_s$ as the maximal pole of order $\geq q+s-1$ for $s \in [1, d]$.\\
The other difference in the argument lies in the fact that in the Stokes formula the  extra term given by the differential of $\psi$ equal to:\\
$$G(\lambda) := \frac{(-1)^{n}}{\Gamma(\lambda)}.\int_U \vert f\vert^{2\lambda}.\bar f^h.u\wedge d\psi $$
has no pole of order $\geq q$ at a point in $\xi + \mathbb{Z}$ because we assumed $d\psi \equiv 0$ near $0$ and we may apply our hypothesis $H(\xi, q)$ with $\varphi := \bar f^h.u\wedge d\psi$ for $h \geq 0$ or with $\varphi := f^{-h}.u\wedge d\psi$ for $h \leq 0$. $\hfill \blacksquare$\\

Let  specialize now the form $\psi \in \mathscr{C}_c^{\infty, (0, n+1)}(U)$ such that $d\psi \equiv 0$ near the origin by defining  $\psi := \rho.\bar \omega'$ where $\omega'$ is a fixed holomorphic $(n+1)-$differential form on $U$ and where $\rho$ is a function in $\mathscr{C}_c^\infty(U)$ which is identically equal to $1$ near the origin  (so $d(\rho.\bar \omega') = d'\rho\wedge \bar \omega'$ vanishes identically near the origin). Then we  consider the Mellin transform of the hermitian period
$$z \mapsto  \frac{1}{(2i\pi)^n}.\int_{f=z} \rho.(\omega/df) \wedge (\overline{\omega'/df}) .$$
In the following corollary we use the hermitian symmetry between $\omega$ and $\omega'$ in order to obtain a better control of the poles of the Mellin transform using the Bernstein polynomials of the frescos associated to $\omega$ and $\omega'$ at the origin.

  \begin{cor}\label{fond.3}
  Let $\tilde{f} : (\mathbb{C}^{n+1}, 0) \to (\mathbb{C},0)$ be a non constant holomorphic germ. Fix $\xi \in \mathbb{Q}$ and a positive integer $q$. Assume that  the hypothesis $H(\xi,q)$ holds for $\tilde{f}$ and consider $\omega$ and $\omega'$ two germs of  $(n+1)-$holomorphic forms. Let $q+d-1, d \geq 1,$ be the maximal order of pole\footnote{Note that the polar parts of order $\geq q$ of the  poles in $\xi + \mathbb{Z}$  are independent of the choice of $\rho$ because of our hypothesis $H(\xi, q)$.} for 
  $$(2i\pi)^{n+1}.F^h_{\omega, \omega'}(\lambda) := \frac{1}{\Gamma(\lambda)}.\int_U \vert f\vert^{2\lambda}.\bar f^h.\rho.\omega\wedge\bar \omega'$$
  for any choice of a point in $\xi + \mathbb{Z}$ and  for any choice of  $h \in \mathbb{Z}$.\\
 Let $\xi_{0}$  be maximal in $\xi + \mathbb{Z}$ such that there exists some $h \in \mathbb{Z}$ and  a pole of order $q+d-1$ at $\xi_{0} \in \xi + \mathbb{Z}$ for $F^h_{\omega, \omega'}(\lambda)$. 
 Then $\xi_0$ is a root of the Bernstein polynomial of the fresco $F_\omega$ and there exists at least $d$ roots of the Bernstein polynomial of $F_\omega$ in $(\xi_0 + \mathbb{N}) \cap [\xi_0, 0[$ counting multiplicities.\\
   Moreover, under our hypothesis, there exists $\xi_1 \in \xi + \mathbb{Z}$ such that for some $h \in \mathbb{Z}$, $F^h_{\omega', \omega}(\lambda)$ has a pole of order $q+d-1$  at \  $\xi_{1}$. Let $\xi_{1} $ be maximal  in $\xi + \mathbb{Z}$ such that this happens. Then  $\xi_1$ is a root of the Bernstein polynomial of the fresco  $F_{\omega' }$ and there exists  at least $d$ roots of the Bernstein polynomial of $F_{\omega' }$ in $(\xi_1 + \mathbb{N}) \cap [\xi_1, 0[$ counting multiplicities.
  \end{cor}
  
  \parag{proof} The first statement is a special case of the previous theorem.\\
  We shall deduce the second statement by using complex conjugaison. Let $\xi_0 \in \xi + \mathbb{Z}$ and $h_0 \in \mathbb{Z}$ such that $F(\lambda) := (2i\pi)^{n+1}.F^{h_0}_{\omega, \omega'}(\lambda)$ has a pole of order $q+d-1$ at $\xi_0$. As $F(\lambda)$ has only real poles, the poles of $\overline{F(\bar \lambda)}$ are the same than the poles of $F(\lambda)$ with the same orders. Moreover we may assume that the function $\rho$ is real, so $\overline{F(\bar \lambda)}$    is  given by
    \begin{align*}
    & \frac{(2i\pi)^{-(n+1)}}{\Gamma(\lambda)}.\int_{U} \vert f\vert^{2(\lambda+h_0)}.\bar f^{-h_0}.\rho.\omega'\wedge\bar \omega =  \frac{(2i\pi)^{-(n+1)}}{\Gamma(\mu-h_0)}.\int_{U} \vert f\vert^{2\mu}.\bar f^{-h_0}.\rho.\omega'\wedge\bar \omega . 
    \end{align*}
    where $\mu = \lambda + h_0$. But $F(\lambda)$ is holomorphic when $\Re(\lambda) \geq 0 $ and $\Re(\lambda + h_0) \geq 0$ so we may replace $\Gamma(\mu -h_0)$ by $\Gamma(\mu)$ in the right hand-side without changing the poles and their orders. We conclude that $F^{-h_0}_{\omega',\omega}$ has a pole of order $q+d-1$ at $ \xi_0 + h_0$ and applying the first statement gives the conclusion.$\hfill\blacksquare$\\

\parag{Remark}

 Of course, with the same method, we can obtain an  analog result than in the theorem \ref{tool 4} for the asymptotic expansion at the origin of a period-integral  of the type
$$ s \mapsto  \int_{\gamma_s} \omega/df $$
where $\omega$ is a holomorphic  $(n+1)-$form on $U$  and where $(\gamma_s)_{s \in H}$ is a horizontal family of compact  $n-$cycles in the fibers of $f$:\\
 Assuming the hypothesis $H(\xi, q)$ for $q \geq 1$ and the existence of a non zero term like $s^{m-\xi}.(Log\, s)^{q+d-2}$ with $d \geq 1$ (or $s^m.(Log\, s)^{q+d-1}$ for $\xi = 0$),  in such an expansion will imply that the Bernstein polynomial of the fresco $F_{\omega}$ will have at least $d$ roots (counting multiplicities) in the set $\xi + \mathbb{Z}$.\\
 This gives a numerical criterion to insure that such a term will not appear in the expansion we are interested in.

\section{The case of a polynomial with $(n+1)$ variables and $(n+2)$ monomials}

The purpose of this section is to give a general algorithm in order to obtain an "estimate" of the Bernstein polynomial of the fresco associated to $(f, \omega)$ for any polynomial $f \in \mathbb{C}[x_0, \dots, x_n]$ with $(n+2)$ monomials  and for  any monomial holomorphic  differential form $\omega = x^\beta.dx$ of degree $n+1$, where $\beta $ is in $\mathbb{N}^{n+1}$ (we exclude the quasi-homogeneous case for $f$ which is obvious). Using the results of the previous sections we obtain a rather precise information of the exponents of the asymptotic expansions of the period integrals $\int_{\gamma_z} \omega/df $ where $(\gamma_z)_{z \in H}$ is a horizontal family of $n-$cycles in the fibers of $f$. This gives also rather precise informations on the poles of the meromorphic extensions of the Mellin transform of the hermitian periods $z \mapsto \int_{f=z} \rho.(\omega/df) \wedge (\overline{\omega'/df})$:
$$ \frac{1}{\Gamma(\lambda)}.\int_{\mathbb{C}^{n+1}} \vert f \vert^{2\lambda}.\bar f^h.\rho.\omega\wedge \bar \omega' $$
where $\rho \in \mathscr{C}^\infty_c(\mathbb{C}^{n+1})$ satisfies $\rho \equiv 1$ near $0$, where $\omega, \omega'$ are monomial holomorphic differential forms of degree $n+1$ and where $h$ is in $\mathbb{Z}$. We shall illustrate the result by several examples.

\subsection{Our setting}

\bigskip

We consider a polynomial \ $f \in \mathbb{C}[x_{0}, \dots, x_{n}] $ \ which is the sum of \ $n+2$ \ monomials 
 $$ f = \sum_{j=1}^{n+2} \ m_{j}$$
   where \ $m_{j} : = \sigma_{j}.x^{\alpha_{j}}$, with \ $\sigma_{j}\in \mathbb{C}^{*}$ \ and \ $\alpha_{j}\in \mathbb{N}^{n+1}$ are not $0$. Define the matrix with \ $(n+1)$ \ lines and \ $(n+2)$ \ columns \ $M = (\alpha_{i,j})$ \ and let  \ $\tilde{M}$ \ be the square  $(n+2, n+2)$ matrix obtained from \ $M$ \ by adding a first line equal to \ $(1, \dots, 1)$.
We shall assume the following conditions:
\begin{enumerate}[(C1)]
\item $\alpha_1, \dots, \alpha_{n+1}$ is a $\mathbb{Q}-$basis of $\mathbb{Q}^{n+1}$.
\item The rank of \ $\tilde{M}$ \ is \ $n+2$.
\end{enumerate}

\parag{Remarks}
\begin{enumerate} 
\item Only the condition (C2) is restrictive on $f$:  when (C2) is fulfilled the condition (C1) may always be satisfied without changing $f$ by a suitable ordering of the $n+2$ monomials.
\item The  condition (C2)  is equivalent to the fact  that \ $f$ \ is {\bf not quasi-homogeneous}.
\end{enumerate}

A diagonal linear change of variables allows to reduce the study to the case where 
\begin{equation*}
 f(x) = \sum_{j=1}^{n+1} \ x^{\alpha_{j}} \ + \lambda.x^{\alpha_{n+2}} \tag{a}
 \end{equation*}
for some \ $\lambda \in \mathbb{C}^{*}$. So in what follows, we shall assume that \ $m_{j} = x^{\alpha_{j}}$ \ for \ $j \in [1,n+1]$ \ and \ $m_{n+2} = \lambda.x^{\alpha_{n+2}}$ \ where \ $\lambda \in \mathbb{C}^{*}$ \ is a parameter.\\

Then we shall write, using our hypothesis (C1):
\begin{equation*}
\alpha_{n+2} = \sum_{j=1}^{n+1} \ \rho_{j}.\alpha_{j}  \tag{b}
\end{equation*}
where $\rho_j$ are rational numbers. We define
\begin{align*}
& H : = \{j \in [1,n+1] \ / \ \rho_{j} = 0 \} ;\\
& J_{+} : = \{j \in [1,n+1] \ / \ \rho_{j} > 0 \} ;\\
& J_{-} : = \{j \in [1,n+1] \ / \ \rho_{j} < 0 \}.
\end{align*}
Let \ $\vert r\vert$ be the smallest positive integer such that \ $\vert r\vert.\rho_{j} : = p_{j}$ \ is an integer for each \ $j \in [1,n+1]$. Write now the relation above as
\begin{equation*} \vert r\vert.\alpha_{n+2} + \sum_{j \in J_{-}} (-p_{j}).\alpha_{j} = \sum_{j \in J_{+}} p_{j}.\alpha_{j}. \tag{c}
\end{equation*}
Now define \ $ d+h$ \ and \ $d$ \ as respectively the supremum and infimum of the two numbers  \ $\vert r\vert + \sum_{j \in J_{-}} \ (-p_{j})$ \ and \ $  \sum_{j \in J_{+}}\  p_{j} $.\\
 Then \ $d$ \ and \ $h$ \ are positive :\\
-the non vanishing of \ $d$ \ is consequence of the fact that \ $\vert r\vert \geq 1$ \ and that at least one\ $p_{j}$ \ is positive.\\
-the non vanishing of \ $h$ \ is consequence of the fact that the equality of these two integers would  imply that the first line in \ $\tilde{M}$ \ satisfies the same linear relation \ $(b)$ \ than all the other lines in \ $\tilde{M}$, contradicting our hypothesis \ (C2).\\
The relation $(e)$ above gives the following equality between the monomials \ $(m_{j})_{j\in [1,n+2]}$:
\begin{equation*}
 m_{n+2}^{\vert r\vert}.\prod_{j \in J_{-}} m_{j}^{-p_{j}} = \lambda^{\vert r\vert}.\prod_{j \in J_{+}} m_{j}^{p_{j}}  \tag{d}
\end{equation*}
and we shall write it 
 \begin{equation*}
  m^{\Delta} = \lambda^{r}.m^{\delta} \tag{e}
  \end{equation*}
 where \ $\Delta$ \ and \ $\delta$ \ are in \ $\mathbb{N}^{n+2}$ \ of respective length \ $d+h$ \ and \ $d$.\\
  Remark that \ $\Delta_j$ \ and \ $\delta_j$ \ are  zero for each \ $j \in H$.\\
  Note that the relation (e) defines the sign of $r$ which is  in \ $\mathbb{Z}^{*}$.\\
  
We shall also use the following observation later on :
\begin{lemma}\label{observation}
The \ $j-$th element of the first column of the matrix \ $\tilde{M}^{-1}$ \ is zero if and only if \ $j$ \ is in \ $H$.
\end{lemma}

\parag{Proof} The co-factor of the element $(1, j)$ in $\tilde{M}$ is the $(n+1, n+1)$ determinant of the matrix with columns $\alpha_1, \dots, \hat{\alpha}_j, \dots, \alpha_{n+2}$. This matrix has rank at most $n$ if and only if $\alpha_{n+2}$ is a linear combination of $\alpha_1, \dots, \hat{\alpha}_j, \dots, \alpha_{n+1}$. This is the case if and only if $\rho_j = 0$, thanks to our hypothesis (C1).$\hfill \blacksquare$\\

\subsection{The result}

Let $\Omega^{p}$ be the $\mathbb{C}[x_0, \dots, x_n]-$module of algebraic $p-$differential forms on $\mathbb{C}^{n+1}$ and fix a polynomial $ f \in \mathbb{C}[x_0, \dots, x_n]$ with $(n+2)$ monomials
\begin{equation}
f := \sum_{j= 1}^{n+2} m_j 
\end{equation}
where $ m_j := x^{\alpha{_j}} $ for $j \in [1, n+1]$ and $m_{n+2} := \lambda.x^{\alpha_{n+1}}$, with $\lambda \in \mathbb{C}^*$, satisfying  the conditions (C1) and (C2) above.\\
 
 We define  $df^p : \Omega^p \to \Omega^{p+1}$ the $\mathbb{C}[x_0, \dots, x_n]-$linear map given by exterior product by $df$ and we note $Ker(df)^p$ its kernel.\\
 Let  $E_f := \Omega^{n+1}\big/ d(Ker(df)^{n}\big)$ endowed with its natural structure of module over the $\mathbb{C}-$algebra $\mathcal{A} := \mathbb{C}\langle a, b\rangle $ where 
 the variables $a$ and $b$ satisfy the commutation relation $ a.b - b.a = b^2$.  Recall that on $E_f$ the action of $a$ is the multiplication by $f$ and the action of $b$ is given by $df\wedge d^{-1}$ where $d$ is the de Rham differential (which is surjective on $\Omega^{n+1}$).\\
 We extend this structure of (left) $\mathcal{A}-$module to $E_f[\lambda] := E_f \otimes_{\mathbb{C}} \mathbb{C}[\lambda]$ by asking that $a$ and $b$ are $\mathbb{C}[\lambda]-$linear.
 
 \begin{lemma}\label{fresco 0}
 For any monomial $x^\beta$  the image of $\mathbb{C}[m_1, \dots, m_{n+2}].x^\beta$ via the map defined by
  $$m^\eta \mapsto \lambda^{\eta_{n+2}}.\prod_{i=0}^{n} x_i^{\sum_{j=1}^{n+2} \alpha_{i, j}.\eta_{j}}  \quad {\rm where} \quad \eta \in \mathbb{N}^{n+2},   $$ 
   is a sub-$\mathcal{A}$-module of $E_f[\lambda]$.
 \end{lemma}
 
 \parag{Proof} We want to prove that this image is stable by the action of $a$ and $b$.\\
 The stability by $a$ is obvious because $a$ is given by the multiplication by $\sum_{j=1}^{n+2} m_j$. \\
  Let $\gamma \in \mathbb{N}^{n+1}$ and compute $b$ using a primitive in $x_i$, for any $ i \in [0, n]$:
 $$b(x^{\gamma+\beta}.dx) = \frac{1}{\gamma_i + \beta_i +1}.x_i.x^{\gamma+\beta}.\frac{\partial f}{\partial x_i}.dx .$$
 Remark now that $x_i.\frac{\partial f}{\partial x_i}.dx = \sum_{j=1}^{n+2} \alpha_{i, j}.m_j $ so the previous computation gives, with $\gamma := M.\eta$ for some $\eta \in \mathbb{N}^{n+2}$ and then  $\gamma_i := \sum_{j=1}^{n+2} \alpha_{i, j}.\eta_j$,
 \begin{equation*}
 \Gamma_i(\eta, \beta).b(m^\eta.x^\beta.dx) = \sum_{j=1}^{n+2} \alpha_{i, j}.m_j.m^\eta.x^{\beta}.dx \quad \forall i \in [0, n], \tag{$@_i$}
 \end{equation*}
 where $\Gamma_i(\eta, \beta) := 1 + \beta_i + \sum_{j=1}^{n+2} \alpha_{i, j}.\eta_j $.\\
 Note that we have also
 \begin{equation*}
 a(m^{\eta}.x^\beta.dx) = \sum_{j=1}^{n+2} m_j.m^\eta.x^\beta.dx \tag{$@_{n+1}$}
 \end{equation*}
 The formulas $(@_i)$ for $i \in [0, n+1]$ are enough to conclude the proof.$\hfill \blacksquare$
 
 \begin{cor}\label{crucial}
Fix $\beta \in \mathbb{N}^{n+1}$. For each $\eta \in \mathbb{N}^{n+2}$ there exists an  element $P_{\beta, \eta}(a, b)$ in  $ \mathcal{A}$, homogeneous of degree $q := \vert \eta\vert := \sum_{j=1}^{n+2} \eta_j$ in (a,b) 
such that:\\
\begin{enumerate}
\item There exists $c(\beta, \eta) \in \mathbb{Q}^*$ such that $P_{\beta, \eta}(a, b)[x^\beta.dx] = c(\beta, \eta).m^\eta.x^\beta.dx $ in $E_f$.
\item Assuming  that $\eta$ satisfies $\eta_j = 0$ for each $j \in H$, there exists rational numbers (depending on $\beta$ and $\eta$) $r_1, \dots, r_q$ such that $P_{\beta, \eta}(a, b) = \prod_{h=1}^q (a - r_h.b)$ in $\mathcal{A}$.
\end{enumerate}
\end{cor}

\parag{Proof} Let first show by induction on $q \geq 0$ that such a $P_{\beta, \eta}(a, b)$ satisfying $1.$ exists and that it satisfies $2.$ when $\eta$ has no component on $H$. As for $q = 0$ the assertion is clear with $P \equiv 1$, assume that our assertion is proved
 for any $\eta$ with $\vert \eta \vert = q-1$ with $q \geq 1$.\\
Then it is enough to prove the assertion for $m_j.m^\eta$ for each $j \in [1, n+2]$ and each $\eta$ with $\vert \eta \vert = q-1$.\\
Then consider the equations $(@_i)$ for $i \in [0, n+1]$ as a square $\mathbb{Q}-$linear system of size $(n+2, n+2)$ with unknown the elements $m_j.m^\eta.x^\beta.dx$ in $E_f$.
The matrix of this system is in $Gl(\mathbb{Q}, n+2)$ thanks to our hypothesis (C2), and so there exists rationals numbers $u_{j}$ and $v_j$ such that we have , for each $j \in [1, n+2]$
$$ m_j.m^\eta.x^\beta.dx = (u_j.a + v_j.b)(m^\eta.x\beta.dx) .$$
With our induction hypothesis this gives that  the homogeneous degree $q$ element in $\mathcal{A}$ defined by $P_{\beta, \eta+1_j}(a, b) := (u_j.a + v_j.b).P_{\beta, \eta}(a,b)$ satisfies $1.$ \\
Assuming that $\eta+1_j$ has no component on $H$ we obtain that $P_{\beta, \eta}(a,b)$ is monic in $a$ (up to a non zero rational number) and then it is enough to show that $u_j$ is not zero to conclude the induction. This is given by the lemma \ref{observation}.$\hfill \blacksquare$\\

\begin{thm}\label{n+2}
Assume that $f \in \mathbb{C}[x_0, \dots, x_n]$ has $(n+2)$ monomials and satisfies the conditions (C1) and (C2) described above. Let $d$ and $h$ the positive integers defined after (c) and $r \in \mathbb{Z}^*$ defined in (c) and (e). For each $\beta \in \mathbb{N}^{n+1}$ there exists homogeneous elements $P_{d+h}$ and $P_d$ of respective degrees $d+h$ and $d$ which are products of homogeneous factors of degree $1$ of the form $a - \xi.b$ where $\xi \in \mathbb{Q}$ such that
\begin{equation}
\Big(P_{d+h}(a, b) - c.\lambda^r.P_d(a, b)\Big)[x^\beta.dx] = 0 \quad {\rm in} \quad E_f[\lambda] \quad {\rm where} \quad c \in \mathbb{Q}^*
\end{equation}
\end{thm}

\parag{Proof}  The previous corollary applied  to both sides of the equality in $E_f[\lambda]$ \\ $m^\Delta.x^\beta.dx = \lambda^r.m^\delta.x^\beta.dx$ deduced from (e) and the corollary \ref{crucial} allow to conclude because we know that $\Delta_j = 0$ and $\delta_j = 0$ for each $j \in H$.$\hfill \blacksquare$\\

This theorem has the following corollary.

\begin{cor}\label{eq. diff.}
Let $f \in \mathbb{C}[x_0, \dots, x_n] $ as in the previous theorem and choose any monomial $x^\beta$.  Define for any horizontal family of $n-$cycles $(\gamma_s)_{s \in S}$ over a  simply connected open set $S$ in $\mathbb{C}^*$ avoiding the critical values of $f$ and having $0$ in its boundary, the integral-period:
\begin{equation}
\varphi_\beta(s) := \int_{\gamma_s}  x^\beta.dx\big/df 
\end{equation}
Then $\varphi_\beta$ is solution on $S$  of the  differential equation (which is regular singular at $0$)  obtained from $(11)$ by letting $a = \times s$ and $b := \int_0^s$. 
\end{cor}

\parag{Proof} Thanks to the proposition \ref{egalite} the Bernstein element of the fresco generated by $[x^\beta.dx]$ in $E_f \otimes_\mathcal{A} \tilde{\mathcal{A}}$ is a left multiple of  $P_d(a, b)$ as $c$ and $\lambda$ are not $0$. This is enough to conclude.$\hfill \blacksquare$\\

This shows that the computation of $P_d(a, b)$ gives rather precise informations on the asymptotic expansion at $0$ of such a  integral-period.\\
We leave the corresponding statement for the poles of the Mellin transform of the hermitian period-integrals corresponding to such an $f$ and monomials holomorphic differential forms $\omega$ and $\omega'$ to the reader. Of course, by conjugaison, the Bernstein polynomial of the fresco $(f, \omega')$ also gives constraints on the possible poles of this Mellin transform, as in the corollary \ref{fond.3}.

 \subsection{Examples}

 The control of the Bernstein polynomial of a fresco will use the theorem \ref{n+2} and the  proposition \ref{egalite}.

 \subsubsection{  $f_{\lambda} := x^{5} + y^{5} + z^{5} + \lambda.x.y.z^{2}$}
 
 We assume that $\lambda$ is a non zero complex number. Then $0$ is the only singular point of the hypersurface $\{ f = 0 \}$ : \\
 as on the set $\Sigma := \{ df = 0\} \subset  \mathbb{C}^{3}$ we have $f(x,y,z) = \frac{1}{5}\lambda.x.y.z $, we easily deduced that $\Sigma \cap \{ f = 0 \} = \{0\}$.
 
 \bigskip
 
 Now using the method explained above we obtain, after some elementary computations, for each monomial form $\omega$ below a degree 4 polynomial multiple of  the Bernstein polynomial of the fresco $F_{\omega}$.\\
 Note that such a fresco has rank at most equal to $4$ and if the rank is equal to $4$ then we obtain the Bernstein polynomial itself.\\
  Of course, the reader interested by more monomials can easily complete this list, where $|$ means ``divides'' :\\
  
 \begin{itemize}
 \item  $\omega = dx\wedge dy\wedge dz \quad \quad  B_{1}(\xi)\, |\, (\xi + \frac{7}{10})(\xi + \frac{4}{5})^{2}(\xi + \frac{6}{5}) $.
 \item  $\omega = x.dx\wedge dy\wedge dz \quad \quad  B_{x}(\xi) \, | \, (\xi + \frac{9}{10})(\xi + 1)(\xi + \frac{6}{5})(\xi + \frac{7}{5}) $.
 \item  $\omega = z.dx\wedge dy\wedge dz \quad \quad  B_{z}(\xi) \, | \, (\xi + 1)^{3}(x + \frac{3}{2}) $.
 \item  $\omega = z^{2}.dx\wedge dy\wedge dz \quad  \quad  B_{z^{2}}(\xi) \, | \, (\xi + \frac{6}{5})^{2} (\xi + \frac{13}{10})(\xi + \frac{9}{5}) $.
  \item  $\omega = x.y.dx\wedge dy\wedge dz \quad \quad  B_{x.y}(\xi) \, | \, (\xi + \frac{11}{10})(\xi + \frac{7}{5})^{2}(\xi + \frac{8}{5})$.
  \item   $\omega = x^{2}.dx\wedge dy\wedge dz \quad \quad  B_{x^{2}}(\xi) \, |\, (\xi + \frac{6}{5})(\xi + \frac{8}{5})^{2} (\xi + \frac{11}{10})$.
  \item  $\omega = x.z.dx\wedge dy\wedge dz \quad \quad  B_{x.z}(\xi) \, | \,(\xi + \frac{6}{5})^{2} (\xi + \frac{7}{5}) (\xi + \frac{17}{10})$.
    \item  $\omega = x.y.z.dx\wedge dy\wedge dz \quad  \quad  B_{x.y.z}(\xi) \, | \, (\xi + \frac{7}{5})(\xi + \frac{8}{5})^{2}(\xi + \frac{19}{10}) \quad$ etc ...
    \end{itemize}
 Note that in this example the differential forms corresponding to degree $2$ monomials in $x, y, z$ are global holomorphic $3-$forms on the fibers of the family of compact surfaces given, for $\lambda$ fixed, by the fibers of the map $\pi_{\lambda}((s,(x, y, z, t)) = s $, sending
 $$ \mathcal{X}_{\lambda} := \{(s,(x, y, z, t)) \in \mathbb{C}\times \mathbb{P}_{3}(\mathbb{C}) \ / \  s.t^{5} = x^{5} + y^{5} + z^{5} +\lambda.x.y.z^{2}.t \} $$
  to $\mathbb{C}$.  As, moreover, the map $\pi_{\lambda}$ has no singular point at infinity, the affine computation controls also the global case for these forms.\\
   Remark that the global computation for these forms gives here the same frescos than in the affine case  because $f_{\lambda}$ has an isolated singularity at the origin.

    \subsubsection{ $f = x.y^{3}+ y.z^{3}+ z.x^{3}+ \lambda.x.y.z$}
    
    The singularity of the hypersurface $\{ f = 0\}$ is the origin :\\
    It is easy to see that any monomial of $f$ is a linear combination of $f$ and $x.\frac{\partial f}{\partial x}, y.\frac{\partial f}{\partial y}, z.\frac{\partial f}{\partial z}$, so that each monomial in $f$ has to vanish on the singular set of $\{ f = 0\}$. Then this implies easily our claim.\\
    Again using the  theorem \ref{n+2} and the  proposition \ref{egalite} allows, after some elementary computations, to find for each monomial form $\omega$ below a degree 3 polynomial dividing the Bernstein polynomial of the fresco $F_{\omega}$.\\
    \begin{itemize}
    \item $\omega = dx\wedge dy\wedge dz \quad \quad  B_{1}(\xi)   \, | \,  (\xi+1)^{3}.$
    \item $\omega = x.dx\wedge dy\wedge dz \quad \quad  B_{x}(\xi)  \, | \,   (\xi + \frac{8}{7})(\xi + \frac{9}{7})(\xi + \frac{11}{7}).$
     \item $\omega = x^{2}.dx\wedge dy\wedge dz \quad \quad  B_{x^{2}}(\xi)   \, | \,   (\xi + \frac{9}{7})(\xi + \frac{11}{7})(\xi + \frac{15}{7}).$
     \item $\omega = x.y.dx\wedge dy\wedge dz \quad \quad  B_{x.y}(\xi)   \, | \,   (\xi + \frac{10}{7})(\xi + \frac{12}{7})(\xi + \frac{13}{7}).$
     \item $\omega = x.y.z.dx\wedge dy\wedge dz \quad \quad  B_{x.y.z}(\xi)   \, | \,  (\xi+2)^{3}.$
      \item $\omega = x^{7}.dx\wedge dy\wedge dz \quad \quad  B_{x^{7}}(\xi)   \, | \,  (\xi+5)(\xi+3).(\xi+2).$ 
      \end{itemize}

      \subsubsection{ $f := x.y^{2}.z^{3} + y.z^{2}.t^{3} + z.t^{2}.x^{3} + t.x^{2}.y^{3} + \lambda.x.y.z.t$}
      
      In this case the singularity is not isolated : the singular of $\{ f = 0 \}$ is the union of the lines $\{ x = y = z = 0 \}, \{y = z = t = 0 \}, \{z = t = x = 0 \}, \{ t = x = y = 0 \}$. The estimate for the  Bernstein polynomial associated to the monomial $1$ (so of the fresco $F_\omega$ with  $\omega := dx\wedge dy\wedge dz\wedge dt$) is $B_{1}(\xi)\, | \, (\xi + 1)^{4}$. So we may have a maximal unipotent monodromy.\\

    \subsubsection{ $f := x.y^{2} + x^{2}.y + z.t^{3} + t.z^{3} + \lambda.x.y.z.t $}
    
    Again we assume that  $\lambda$ is a non zero complex number. The hypersurface $ \{ f = 0 \} $ has an isolated singularity at the origin :\\
    If $\Sigma := \{ df = 0 \} \subset \mathbb{C}^{4}$ we have on $\Sigma$ the relations $x.y^{2} = x^{2}.y = \frac{-1}{3}.\lambda.x.y.z.t$ and $z.t^{3}= z^{3}.t =  \frac{-1}{4}.\lambda.x.y.z.t$. So on $\Sigma \cap \{f = 0\}$ we have $x.y = 0 = z.t $ and this implies that $\Sigma \cap \{f = 0\} = \{0\}$.\\
    
    Now we shall use again the  theorem \ref{n+2} and the  proposition \ref{egalite} in order to give a polynomial of degree $12$ which divides the Bernstein polynomial of the fresco $F_{\omega}$ for $\omega := dx \wedge dy\wedge dz\wedge dt$. The reader interested by another holomorphic monomial form can follow the same line to obtain an analogous result.\\
    
    The relation between the monomials of $f$ is 
    $$\lambda^{12}(x.y^{2})^{4}(y.x^{2})^{4}(z.t^{3})^{3}(z^{3})^{3} = (\lambda.x.y.z.t)^{12}.$$
    So to compute the initial form in (a,b) of the polynomial in $\mathcal{A}$ constructed in the theorem \ref{n+2} annihilating $[\omega]$ in $E^{4}\big/B(E^{4})$, it is enough to compute the homogeneous in (a,b) polynomial $P$ of degree $12$ satisfying in $E^{4}$ the relation  $P.[\omega] = [(\lambda.x.y.z.t)^{12}.\omega]$. \\
    Note $m_{1}, \dots, m_{4}$ the first monomials in $f$ and $m := \lambda.x.y.z.t $. Then we have in $E^{4}$ the equality for any integer $k \geq 0$ (where $\omega$ is omitted)\\
    \begin{itemize}
    \item $m_{1}.m^{k} =  \frac{1}{3}.\big((k+1).b[m^{k}] - m^{k+1}] $ 
    \item  $m_{2}.m^{k}  = \frac{1}{3}.\big((k+1).b[m^{k}] - m^{k+1}] $ 
    \item $ m_{3}.m^{k}  = \frac{1}{4}.\big((k+1).b[m^{k}] - m^{k+1}] $
     \item $ m_{4}.m^{k}  = \frac{1}{4}.\big((k+1).b[m^{k}] - m^{k+1}] $
     \end{itemize}
     and so we obtain
     $$ \big(a - \frac{7}{6 }(k+1).b\big)[m^{k}] = \frac{-1}{6}m^{k+1} .$$
     Then the initial form of the polynomial annihilating $[\omega]$ is equal to the product ordered from left to right by decreasing $k$
     $$ \prod_{k=0}^{11} \ \big(a - \frac{7}{6}(k+1).b\big)[m^{k}]  .$$
     This gives the following estimate for the Bernstein polynomial
     $$ B(\xi)  \, | \,   \prod_{k=0}^{11} \ (\xi +  \frac{k+7}{6})    $$


\begin{thebibliography}{99}
 
 \parag{References}
 
   \bibitem{[B.81]} Barlet, D. {\it  D\'eveloppements asymptotiques des fonctions obtenues par int\'egration sur les fibres},  Inv. Math. vol. 68 (1982), p. 129-174.
  \bibitem{[B-M.87]} Barlet, D. et Maire, H.M. \textit{D\'eveloppements asymptotiques, transformation de Mellin complexe et int\'egration dans les fibres},  in Sem. P. Lelong, Lecture Notes, vol. 1295 Springer Verlag  (1987), p. 11-23.
  \bibitem{[B.91]}  Barlet, D. {\it Interaction de strates cons\'ecutives pour les cycles \'evanescents}, Ann. Sci. \'Ecole Norm. Sup. (4) 24 (1991), no. 4, pp. 401-505.
  \bibitem{[B.93]} Barlet, D. {\it Th\'eorie des (a,b)-modules I}, Complex Analysis and Geometry, pp. 1-43, Plenum Press, New York 1993.
 \bibitem{[B.06]} Barlet, D. \textit{ Sur certaines singularit\'es non isol\'ees d'hypersurfaces I}, Bull. Soc. Math. France 134 fasc.2  (2006), p. 173-200.
  \bibitem{[B-S.07]} Barlet, D. and Saito, M. \textit{Brieskorn modules and Gauss-Manin systems for non isolated hypersurface singularities} Bull. of London Math. Soc. (2007).
  \bibitem{[B.08]} Barlet, D. \textit{Sur certaines singularit\'es d'hypersurfaces II},   Journal of Algebraic  Geometry 17 (2008), p.199-254.
   \bibitem{[B.09a]} Barlet, D. {\it Sur les fonctions \`a lieu singulier de dimension 1}, Bull. Soc. math. France 137 (4), (2009), p. 587-612.
  \bibitem{[B.09b]} Barlet, D. {\it P\'eriodes \'evanescentes et (a,b)-modules monog\`enes}, Bollettino U.M.I. (9) II (2009) p. 651-697.
  \bibitem{[B.13a]} Barlet, D. {\it Asymptotics of a vanishing period : characterization of semi-simplicity},  arXiv:1301.7589,  math.AG and math.CV.
  \bibitem{[B.12]} Barlet, D. {\it A finiteness theorem for S-relative formal Brieskorn module}, math. arXiv 1207.4013, math.AG and math.CV.
  \bibitem{[B.13]} Barlet, D. {\it  Algebraic differential equations associated to some polynomials}, math. arXiv:1305.6778, math.AG and math.CV. 
  \bibitem{[B.15]} Barlet, D. {\it A note on some fiber-integrals} math. arXiv:1512.07062, math.CV and  math.AG:
  \bibitem{[Br.70]} Brieskorn, E. {\it Die Monodromie der isolierten Singularit\"aten von Hyperfl\"achen}, Manuscripta Math. 2 (1970), pp. 103-161.
 \bibitem{[De.70]} Deligne, P. {\it Equations diff\'erentielles \`{a} points singuliers r\'eguliers}, L-N 163 (1970) Springer.
  \bibitem{[K.76]} Kashiwara, M. {\it b-function and holonomic systems, rationality of roots of b-functions}, Invent. Math. 38 (1976) p. 33-53.
  \bibitem{[Ma.74]} Malgrange, B. {Int\'egrale asymptotique et monodromie}, Ann. Sc. ENS  t.7 (1974), pp.405-430.
  \bibitem{[Mi.68]} Milnor, J. {\it Singular Points of Complex Hypersurfaces}, Ann. of Math. Studies 61 (1968) Princeton.
  \bibitem{[Se.70]} Sebastiani, M. {\it Preuve d'une conjecture de Brieskorn}, Manuscripta Math. 2 (1970), pp. 301-308.
  \bibitem{[Sk.83]} Saito, Kyoji {\it Period mapping associated to a primitive form}, Publ. Res. Inst. Math. Sci. 19 (1983), no. 3, pp.1231-1264
  \bibitem{[Sa.89]} Saito, Morihiko {\it On the structure of Brieskorn lattice},  Ann. Inst. Fourier (Grenoble) 39 (1989), no. 1, pp. 27-72.
  \bibitem {[V.80]} Varchenko, A. N. {\it Asymptotic behavior of holomorphic forms determines a mixed Hodge structure}, (Russian) Dokl. Akad. Nauk SSSR 255 (1980), no. 5, pp. 1035-1038. 
 
 \end{thebibliography}
 \end{document}